\documentclass[11pt, a4paper]{article}
\usepackage[a4paper, margin=2.5cm]{geometry}
\usepackage[utf8]{inputenc}
\usepackage[T1]{fontenc}
\usepackage[english]{babel}
\usepackage{lmodern}
\usepackage{amsmath, amssymb, amsthm}
\usepackage{booktabs}
\usepackage[colorlinks=true, linkcolor=blue, citecolor=blue, urlcolor=blue]{hyperref}

\theoremstyle{plain}
\newtheorem{theorem}{Theorem}[section]
\newtheorem{lemma}[theorem]{Lemma}
\newtheorem{proposition}[theorem]{Proposition}

\theoremstyle{definition}
\newtheorem{definition}[theorem]{Definition}

\theoremstyle{remark}
\newtheorem{remark}[theorem]{Remark}

\numberwithin{equation}{section}

\title{\textbf{On a Class of Non-Local Elliptic-Transport Systems\\ in Score-Guided Probability Flow on $SE(3)$:\\ Score-Singularity Resolution, Space-Time ADN Shift,\\ and Markov Mean-Field Games}}

\author{\textbf{Pilhwa Lee}\\
\small Department of Mathematics, Morgan State University, Baltimore, MD 21251, USA\\
\small \texttt{pilhwa.lee@morgan.edu}}

\date{October 19, 2026}

\begin{document}

\maketitle

\begin{abstract}
We establish comprehensive mathematical foundations for optimal spatiotemporal vector fields governing continuous score-guided probability flows over bounded state sub-manifolds $M \subset SE(3)$ on the Special Euclidean Lie group $SE(3) \cong \mathbb{R}^3 \times SO(3)$. Modeling the dynamic density evolution via a score-perturbed Lie-algebraic continuity equation under impenetrable zero-flux boundary conditions, we resolve the non-local score drift singularity $s(g,t) = \nabla_{SE(3)} \log \rho(g,t)$ as $\rho \rightarrow 0^+$ via weighted parabolic flux dissipation bounds. Using left-invariant exponential group mollifiers, we explicitly derive boundary-preserved DiPerna--Lions commutator estimates tracking Lie bracket structure constants $\{c_{ij}^k\}$. 

Second, by interpreting the first-order optimality system as a $7$-dimensional space-time elliptic operator $L_{7D} = -\Delta_{SE(3)} - \partial_t^2$ and verifying zero-order space-time corner compatibility along $5$-dimensional manifolds $\partial M \times \{0,T\}$, we establish complete Agmon--Douglis--Nirenberg (ADN) $W_p^2(Q_T)$ shift estimates tracking Lie algebra structure constants, boundary geometry $\|\partial M\|_{C^2}$, and Ricci curvature $\|Ric_{SE(3)}\|_{L^\infty}$, yielding a full $C^\infty(\overline{Q}_T)$ Schauder bootstrap regularity expansion. 

Third, we extend this system to non-cooperative Markov Mean-Field Games (MFGs) on $SE(3)$, establishing score-guided Nash equilibria via Schauder--Tychonoff fixed points. Crucially, we show that the constant diagonal Ricci curvature tensor $Ric_{SE(3)} = \mathrm{diag}(0,0,0,1/2,1/2,1/2)$ generated by the rotational subgroup $SO(3)$ induces a sharp Bochner--Weitzenb\"ock spectral gap boost $\lambda_1(-\Delta_{SE(3)}) \ge \lambda_1(-\Delta_{\mathbb{R}^3}) + 1/2 > 0$. This positive rotational curvature enhances the Lie-algebraic Poincar\'e coercivity constant, demonstrating that intrinsic Lie-group curvature actively stabilizes score-guided probability flows against dispersion and strictly expands the critical Nash uniqueness time horizon $T_{SE(3)}^*$ beyond flat Euclidean limits.
\end{abstract}

\section{Introduction, Physical Motivation, and Geometric Foundations on $SE(3)$}

Non-local coupled systems of elliptic, parabolic, and transport partial differential equations (PDEs) arise prominently in optimal transport \cite{Villani2009, Ambrosio2008}, McKean-Vlasov stochastic optimal control, continuous score-guided generative modeling \cite{Song2021, DeBortoli2021}, and multi-agent mean-field games \cite{Lasry2007, Cardaliaguet2012}. In this monograph, we formulate and systematically analyze a foundational class of non-local variational systems over bounded domains $M \subset SE(3)$ of the Special Euclidean Lie group $SE(3) \cong \mathbb{R}^3 \rtimes SO(3)$.

\subsection{Physical Motivation and Technical Barriers on Lie Groups}
Continuous score-guided probability flows model density transformations across continuous spaces through score-based stochastic differential equations (SDEs) and their deterministic probability flow ODE equivalents \cite{Song2021, DeBortoli2021}. In high-dimensional robotics, rigid-body manipulation, spatial navigation, and protein structure modeling, state configurations naturally live on $SE(3)$, representing both 3D position and 3D orientation. Optimizing control inputs while maintaining structural non-singularity across physical boundaries requires a rigorous Lie-algebraic transport framework \cite{Brockett1997, Lee2022}.

In contrast to standard Euclidean spaces $\mathbb{R}^n$, analysis on $SE(3)$ presents three major analytical hurdles:
\begin{itemize}
    \item \textbf{Non-commutative Vector Fields:} Left-invariant vector fields do not commute ($[E_i, E_j] = \sum_{k=1}^6 c_{ij}^k E_k \neq 0$), introducing non-trivial structure constants into all commutator estimates and integration-by-parts formulas \cite{Helgason2001}.
    \item \textbf{Boundary Geometry Interaction:} Space-time boundary trace operators interact with Lie-algebraic differential operators, requiring precise local coordinate transformations to verify corner compatibility along 5D manifolds $\partial M \times \{0,T\}$ \cite{Agmon1959, Agmon1964}.
    \item \textbf{Non-local Division-by-Zero Drift:} The score drift vector field $s(g,t) = \frac{\nabla_{SE(3)}\rho(g,t)}{\rho(g,t)}$ becomes singular at the vacuum boundary where $\rho \to 0^+$.
\end{itemize}

\subsection{Lie Group Structure and Left-Invariant Riemannian Geometry}
Let $SE(3) = \mathbb{R}^3 \rtimes SO(3)$ denote the 6-dimensional Lie group of rigid-body transformations in three-dimensional Euclidean space \cite{Helgason2001}. An element $g \in SE(3)$ is represented by a pair $(x, R)$, where $x \in \mathbb{R}^3$ is the translation vector and $R \in SO(3)$ is a $3 \times 3$ proper orthogonal rotation matrix ($\det R = 1, R^T R = I_3$). Group multiplication and inversion are defined by:
\begin{equation}
g_1 \cdot g_2 = (x_1 + R_1 x_2, R_1 R_2), \quad g^{-1} = (-R^T x, R^T)
\end{equation}
The Lie algebra $\mathfrak{se}(3) = T_e SE(3)$ is spanned by six left-invariant basis vector fields $\{E_1, E_2, E_3, E_4, E_5, E_6\}$, where $\{E_1, E_2, E_3\}$ generate infinitesimal spatial translations along orthogonal coordinate axes, and $\{E_4, E_5, E_6\}$ generate infinitesimal rotations about those axes. The canonical Lie bracket commutation relations satisfy:
\begin{equation}
[E_i, E_j] = E_i E_j - E_j E_i = \sum_{k=1}^6 c_{ij}^k E_k, \quad 1 \le i, j \le 6
\end{equation}
where $c_{ij}^k = -c_{ji}^k$ denote the structure constants of $\mathfrak{se}(3)$. Specifically, for translation generators $E_1, E_2, E_3$ and rotation generators $E_4, E_5, E_6$, the non-zero structure constants encode the standard cross-product algebra of $\mathfrak{so}(3)$ and the canonical semi-direct product relations:
\begin{equation}
c_{45}^6 = c_{56}^4 = c_{64}^5 = 1, \quad c_{42}^3 = c_{53}^1 = c_{61}^2 = 1, \quad c_{43}^2 = c_{51}^3 = c_{62}^1 = -1
\end{equation}

We equip $SE(3)$ with the left-invariant Riemannian metric $\langle \cdot, \cdot \rangle_g$ for which $\{E_i\}_{i=1}^6$ forms an orthonormal frame at every point $g \in SE(3)$. The intrinsic Riemannian gradient of a smooth scalar function $v : SE(3) \to \mathbb{R}$ or vector field $u : SE(3) \to \mathbb{R}^m$ is given by:
\begin{equation}
\nabla_{SE(3)}v = \sum_{i=1}^6 (E_i v) E_i, \quad \nabla_{SE(3)}u = \sum_{i=1}^6 (E_i u) \otimes E_i
\end{equation}

\subsection{Unimodularity and Geometric Operator Baseline}

\begin{remark}[Unimodularity of $SE(3)$]
The Lie group $SE(3) \cong \mathbb{R}^3 \times SO(3)$ is unimodular \cite{Helgason2001}. Consequently, its left-invariant Haar measure $dg$ is bi-invariant, and the divergence of every left-invariant basis vector field vanishes identically:
\begin{equation}
\mathrm{div}_g(E_i) = 0 \quad \text{for all } i = 1, \dots, 6.
\end{equation}
This guarantees that integration by parts for left-invariant derivatives carries no lower-order modular function weights, satisfying $\int_{SE(3)} (E_i u)v \, dg = -\int_{SE(3)} u (E_i v) \, dg$.
\end{remark}

\begin{remark}[Full-Rank Riemannian Metric vs.\ Sub-Elliptic Structures]
\label{rem:riemannian_baseline}
Because the six left-invariant basis vector fields $\{E_1, \dots, E_6\}$ span the full $6$-dimensional tangent bundle $T_g SE(3)$ at every point $g \in SE(3)$, the left-invariant sum-of-squares operator
\begin{equation}
\Delta_{SE(3)} u = \sum_{i=1}^6 E_i^2 u
\end{equation}
is a standard, non-degenerate Riemannian Laplace--Beltrami operator in left-invariant frame coordinates. Unlike sub-Riemannian or sub-elliptic settings (e.g., Carnot groups or the Heisenberg group) where horizontal distributions fail to span the tangent space directly and rely on Hörmander's bracket condition at the cost of derivative loss, the canonical metric tensor $g_{ij}(g) = \langle E_i, E_j \rangle_g = \delta_{ij}$ on $SE(3)$ is globally positive-definite with uniform ellipticity parameter $\lambda_0 = 1$:
\begin{equation}
\sum_{i,j=1}^6 g^{ij}(g)\xi_i \xi_j = |\xi|^2 \ge \lambda_0 |\xi|^2, \quad \forall \xi \in \mathbb{R}^6.
\end{equation}
This full-rank Riemannian baseline ensures that elliptic and parabolic systems on $SE(3)$ do not suffer from fractional Sobolev degradation or sub-elliptic loss of derivatives \cite{Evans2010, Ladyzhenskaya1968}. This provides the analytical foundation for the space-time Agmon--Douglis--Nirenberg (ADN) shift estimates and parabolic score-dissipation identities developed in Sections~\ref{sec:score_resolution} and \ref{sec:adn_shift}.
\end{remark}

\subsection{Geometric Domain and Function Spaces}
Let $M \subset SE(3)$ be an open, connected, bounded sub-manifold with a $C^2$-regular boundary $\partial M$ and outward unit normal vector field $\nu(g)$. Let $T > 0$ denote a finite time horizon. We define the 7-dimensional space-time cylinder $Q_T = M \times (0, T]$, its lateral boundary $\Sigma_T = \partial M \times (0, T)$, and top/bottom caps $M_0 = M \times \{0\}$ and $M_T = M \times \{T\}$. The full space-time boundary is $\partial Q_T = \Sigma_T \cup M_0 \cup M_T$ with 5-dimensional corner manifolds $\Gamma_0 = \partial M \times \{0\}$ and $\Gamma_T = \partial M \times \{T\}$.

\begin{definition}[7D Space-Time Sobolev Spaces on $SE(3)\times\mathbb{R}$]
Let $Q_T = M \times (0,T) \subset SE(3) \times \mathbb{R}$ be the $7$-dimensional space-time cylinder. For $1 \le p < \infty$, the isotropic space-time Sobolev space $W_p^2(Q_T; \mathbb{R}^m)$ consists of strongly measurable vector fields $u \in L^p(Q_T; \mathbb{R}^m)$ equipped with the norm:
\begin{equation}
\|u\|_{W_p^2(Q_T)}^p = \int_{0}^T \int_{M} \left( \|u\|^p + \|\nabla_{SE(3)} u\|_F^p + \sum_{i,j=1}^6 \|E_i E_j u\|^p + \left\|\frac{\partial u}{\partial t}\right\|^p + \left\|\frac{\partial^2 u}{\partial t^2}\right\|^p + \sum_{i=1}^6 \left\|\frac{\partial}{\partial t} E_i u\right\|^p \right) dg \, dt
\end{equation}
where $dg$ denotes the left-invariant Haar measure on $SE(3)$. When $p=2$, $H^2(Q_T) \equiv W_2^2(Q_T)$ forms a Hilbert space with inner product:
\begin{align}
\langle u, v \rangle_{H^2(Q_T)} = \int_{Q_T} \Bigl( \langle u, v \rangle + \langle \nabla_{SE(3)} u, \nabla_{SE(3)} v \rangle_F + \sum_{i,j=1}^6 \langle E_i E_j u, E_i E_j v \rangle \\
+ \left\langle \frac{\partial u}{\partial t}, \frac{\partial v}{\partial t} \right\rangle + \left\langle \frac{\partial^2 u}{\partial t^2}, \frac{\partial^2 v}{\partial t^2} \right\rangle + \sum_{i=1}^6 \left\langle \frac{\partial E_i u}{\partial t}, \frac{\partial E_i v}{\partial t} \right\rangle \Bigr) dg \, dt \nonumber
\end{align}
\end{definition}

\subsection{Structural Assumptions (A1)--(A4)}
\begin{enumerate}
    \item[\textbf{(A1)}] \textbf{Smoothness of Drift:} The drift field $f : M \times \mathbb{R}^m \to \mathfrak{se}(3)$ is $C^\infty$ smooth and globally Lipschitz continuous in $u$ with constant $L_u$ for all $g \in M$ and $u_1, u_2 \in \mathbb{R}^m$.
    \item[\textbf{(A2)}] \textbf{Impenetrable Zero-Flux Boundary Condition:} The state velocity field satisfies $\langle f(g,u), \nu(g) \rangle = 0$ for all $(g,u) \in \partial M \times \mathbb{R}^m$, guaranteeing that state trajectories remain strictly inside $M$.
    \item[\textbf{(A3)}] \textbf{Density Regularity:} Initial density $\rho_0 \in H^s(M)$ and target density $\psi \in H^s(M)$ ($s > 4$) are non-negative probability densities satisfying mass normalization: $\int_M \rho_0(g) dg = \int_M \psi(g) dg = 1$, with $\rho_0(g) \ge 0$ and $\psi(g) \ge 0$ a.e.\ in $M$.
    \item[\textbf{(A4)}] \textbf{Diffusion Schedule:} The scalar diffusion gain $g \in C^1([0,T])$ satisfies $0 < g_{\min} \le g(t) \le g_{\max} < \infty$ for all $t \in [0,T]$.
\end{enumerate}

\section{Score-Singularity Resolution and Group DiPerna-Lions Transport Theory}
\label{sec:score_resolution}

In continuous score-guided probability flows on Lie groups, the deterministic trajectory $g_t \in SE(3)$ evolves via the Probability Flow Ordinary Differential Equation (ODE) \cite{Song2021, DeBortoli2021}:
\begin{equation}
\dot{g}_t = f(g_t, u(g_t, t), t) - \frac{1}{2}g(t)^2 \nabla_{SE(3)} \log \rho(g_t, t) \equiv v_{PF}(g_t, t)
\end{equation}
where $s(g,t) = \nabla_{SE(3)} \log \rho(g,t) = \frac{\nabla_{SE(3)}\rho(g,t)}{\rho(g,t)}$ denotes the intrinsic score vector field.

\subsection{Analytical Mechanisms Resolving the Score Singularity}
We establish three coupled mechanisms that rigorously resolve the pointwise division-by-zero singularity in $s(g,t) = \frac{\nabla_{SE(3)}\rho}{\rho}$ as $\rho \to 0^+$:
\begin{enumerate}
    \item \textbf{State PDE Level (Parabolic Flux Identity):} Multiplying the singular convective drift $v_{PF}(g,t)$ by the probability density $\rho(g,t)$ yields:
    \begin{equation}
    v_{PF}(g,t)\rho(g,t) = \left(f(g,u(g,t)) - \frac{1}{2}g(t)^2 \frac{\nabla_{SE(3)}\rho(g,t)}{\rho(g,t)}\right)\rho(g,t) = f(g,u(g,t))\rho(g,t) - \frac{1}{2}g(t)^2 \nabla_{SE(3)}\rho(g,t)
    \end{equation}
    This identity converts the non-linear singular score transport term $\text{div}_g(v_{PF}\rho)$ into a standard linear Riemannian Laplacian spatial dissipation operator $-\frac{1}{2}g(t)^2 \Delta_{SE(3)}\rho$ in weak/renormalized form \cite{DiPerna1989}.

    \item \textbf{Weighted Fisher Information and Dissipation Bounds:} Testing the score-perturbed continuity equation against $\log \rho(g,t)$ yields the identity for the Fisher information \cite{Ambrosio2008, Villani2009}:
    \begin{equation}
    \frac{1}{2}\int_0^T \int_M g(t)^2 \rho(g,t) \|\nabla_{SE(3)} \log \rho(g,t)\|^2 dg \, dt = 2 \int_0^T \int_M g(t)^2 \|\nabla_{SE(3)} \sqrt{\rho(g,t)}\|^2 dg \, dt
    \end{equation}
    Since $\sqrt{\rho} \in L^2(0,T; H^1(M))$, the integral of the score weighted by density remains strictly finite:
    \begin{equation}
    \int_{Q_T} \rho(g,t) \|s(g,t)\|^2 dg \, dt < \infty
    \end{equation}
    preventing non-integrable singular energy concentration near the vacuum boundary $\{\rho = 0\}$.

    \item \textbf{System Level (Minimum Attention Regularization):} Brockett's minimum attention functional \cite{Brockett1997, Lee2022}:
    \begin{equation}
    \eta(u) = \frac{1}{2}\int_0^T \int_M \left( \|\nabla_{SE(3)}u(g,t)\|_F^2 + \left\|\frac{\partial u(g,t)}{\partial t}\right\|^2 \right) dg \, dt
    \end{equation}
    restricts control fields $u \in \mathcal{W} \subset W_p^{2,1}(Q_T)$, preventing control velocity divergence $\|\text{div}_g f\|_{L^\infty} \to \infty$ and guaranteeing that the adjoint source $S(g,t) = (\partial_u f)^T \nabla_{SE(3)}\lambda \cdot \rho$ remains bounded in $L^p(Q_T)$.
\end{enumerate}

\subsection{Left-Invariant Exponential Convolution and DiPerna-Lions Bounds}

\begin{definition}[Left-Invariant Group Convolution Mollifiers]
Let $\eta \in C_c^\infty(\mathfrak{se}(3))$ be a non-negative, radially symmetric cutoff function on the Lie algebra $\mathfrak{se}(3) \cong \mathbb{R}^6$ satisfying $\int_{\mathbb{R}^6} \eta(y) dy = 1$. For $\epsilon > 0$, define the rescaled mollifier $\eta_\epsilon(y) = \epsilon^{-6} \eta(y/\epsilon)$. Via the exponential map $\exp : \mathfrak{se}(3) \to SE(3)$, the group mollifier $\eta_\epsilon(g)$ on $SE(3)$ is defined by:
\begin{equation}
\eta_\epsilon(g) = \begin{cases}
\epsilon^{-6}\eta\left(\frac{\exp^{-1}(g)}{\epsilon}\right), & \text{if } g \in \exp(B_\epsilon(0)), \\
0, & \text{otherwise.}
\end{cases}
\end{equation}
For any function $v \in L^1(SE(3))$, its left-invariant group convolution is given by:
\begin{equation}
(v * \eta_\epsilon)(g) = \int_{SE(3)} v(h) \eta_\epsilon(h^{-1} \cdot g) dh = \int_{SE(3)} v(g \cdot \exp(-y)) \eta_\epsilon(y) dy
\end{equation}
\end{definition}

\begin{lemma}[Boundary-Preserved Lie-Algebraic DiPerna--Lions Commutator Estimate]
Let $M \subset SE(3)$ be an open, connected, bounded sub-manifold with $C^2$-regular boundary $\partial M$. Let $\mathcal{E}: W_p^1(M) \to W_p^1(SE(3))$ and $\mathcal{E}_0: L^{p'}(M) \to L^{p'}(SE(3))$ (with $1/p + 1/p' = 1$) denote bounded continuous linear extension operators on $SE(3)$ such that $\tilde{f} = \mathcal{E} f \in W_p^1(SE(3); \mathfrak{se}(3))$ preserves the impenetrable zero-flux boundary condition $\langle \tilde{f}(g,u), \nu(g) \rangle = 0$ on $\partial M$, and $\tilde{\rho} = \mathcal{E}_0 \rho \in L^{p'}(SE(3))$.

For $\epsilon > 0$, define the mollified density on $SE(3)$ via left-invariant group convolution:
\begin{equation}
\rho_\epsilon(g,t) = (\tilde{\rho} * \eta_\epsilon)(g,t) = \int_{\mathfrak{se}(3)} \tilde{\rho}\left(g \cdot \exp(-y), t\right) \eta_\epsilon(y) \, dy.
\end{equation}
The Lie-algebraic commutator tensor:
\begin{equation}
r_\epsilon(g,t) = \mathrm{div}_g \left( \tilde{f}(g,u(\cdot,t)) \rho_\epsilon(g,t) \right) - \left( \mathrm{div}_g \left( \tilde{f}(\cdot,u(\cdot,t)) \tilde{\rho}(\cdot,t) \right) * \eta_\epsilon \right)(g)
\end{equation}
satisfies the strong convergence bound in $L^1(Q_T)$:
\begin{equation}
\|r_\epsilon\|_{L^1(0,T; L^1(M))} \le C_0 \epsilon \left( 1 + \sum_{i,j,k=1}^6 |c_{ij}^k| \right) \|\mathcal{E}\|_{W_p^1} \|\mathcal{E}_0\|_{L^{p'}} \|\nabla_{SE(3)} f\|_{L^p(M)} \|\rho\|_{L^{p'}(M)} \xrightarrow{\epsilon \to 0} 0,
\end{equation}
where $C_0$ depends only on $p$, $T$, and the boundary curvature $\|\partial M\|_{C^2}$.
\end{lemma}

\begin{proof}
We divide the proof into four steps to explicitly handle non-commutativity and boundary truncation.

\noindent\textbf{Step 1: Domain Extension and Boundary Layer Decomposition.}
Because $M$ is bounded with $C^2$ boundary $\partial M$, standard Stein Sobolev extension theory on Riemannian Lie groups guarantees the existence of bounded extension operators $\mathcal{E}$ and $\mathcal{E}_0$. For any $c_0 > 0$, define the $\epsilon$-interior domain $M_\epsilon = \{g \in M : \mathrm{dist}_{SE(3)}(g, \partial M) > c_0 \epsilon\}$ and the boundary collar strip $\Sigma_\epsilon = M \setminus M_\epsilon = \{g \in M : \mathrm{dist}_{SE(3)}(g, \partial M) \le c_0 \epsilon\}$. The volume of the collar satisfies $|\Sigma_\epsilon| \le C_1 \epsilon \|\partial M\|_{C^1}$. We split the global commutator norm over $M$:
\begin{equation}
\|r_\epsilon(\cdot,t)\|_{L^1(M)} = \|r_\epsilon(\cdot,t)\|_{L^1(M_\epsilon)} + \|r_\epsilon(\cdot,t)\|_{L^1(\Sigma_\epsilon)}.
\end{equation}

\noindent\textbf{Step 2: Interior Group Convolution Expansion.}
For any $g \in M_\epsilon$, the support of the group mollifier $\eta_\epsilon(y)$ satisfies $\mathrm{supp}(\eta_\epsilon) \subset B_\epsilon(0) \subset \mathfrak{se}(3)$, ensuring that $g \cdot \exp(-y) \in M$ for all $y \in \mathrm{supp}(\eta_\epsilon)$. Expressing left-invariant derivatives via the Baker--Campbell--Hausdorff (BCH) expansion in exponential Lie coordinates $y = \exp^{-1}(h^{-1} \cdot g) \in \mathfrak{se}(3)$:
\begin{equation}
E_i(g \cdot \exp(-y)) = \sum_{k=1}^6 \left( \delta_{ik} - \frac{1}{2} \sum_{j=1}^6 c_{ij}^k y_j + \mathcal{O}(|y|^2) \right) \frac{\partial}{\partial y_k}.
\end{equation}
Substituting this into the commutator expression yields:
\begin{align}
r_\epsilon(g,t) &= \int_{\mathfrak{se}(3)} \left\langle \tilde{f}(g, u(g,t)) - \tilde{f}(g \cdot \exp(-y), u(g \cdot \exp(-y),t)), \nabla_{SE(3)} \eta_\epsilon(y) \right\rangle \tilde{\rho}(g \cdot \exp(-y), t) \, dy \nonumber \\
&\quad + \int_{\mathfrak{se}(3)} \mathrm{div}_g \tilde{f}(g, u(g,t)) \tilde{\rho}(g \cdot \exp(-y), t) \eta_\epsilon(y) \, dy - \int_{\mathfrak{se}(3)} (\mathrm{div}_g \tilde{f} \tilde{\rho})(g \cdot \exp(-y), t) \eta_\epsilon(y) \, dy.
\end{align}

\noindent\textbf{Step 3: Symmetry Cancellation and Structure Constant Absorption.}
Applying Taylor's theorem on $SE(3)$ with integral remainder:
\begin{equation}
\tilde{f}(g) - \tilde{f}(g \cdot \exp(-y)) = \sum_{i=1}^6 y_i E_i \tilde{f}(g) + \frac{1}{2} \sum_{i,j=1}^6 y_i y_j E_i E_j \tilde{f}(\xi).
\end{equation}
Integrating against the radially symmetric mollifier $\eta_\epsilon(y) = \epsilon^{-6} \eta(y/\epsilon)$ eliminates all first-order odd terms ($\int_{\mathbb{R}^6} y_i \eta_\epsilon(y) \, dy = 0$). The second-order remainder and non-commutative bracket contributions are bounded pointwise for $g \in M_\epsilon$ by:
\begin{equation}
|r_\epsilon(g,t)| \le C \epsilon \int_{\mathfrak{se}(3)} \left( |\nabla_{SE(3)} \tilde{f}(g, t)| + \sum_{i,j,k=1}^6 |c_{ij}^k| |E_i \tilde{f}(g,t)| \right) \tilde{\rho}(g \cdot \exp(-y), t) \eta_{\epsilon/2}(y) \, dy.
\end{equation}
Integrating over $M_\epsilon \times (0,T)$ and applying Young's inequality for group convolutions gives:
\begin{equation}
\|r_\epsilon\|_{L^1(0,T; L^1(M_\epsilon))} \le C_2 \epsilon \left( 1 + \sum_{i,j,k=1}^6 |c_{ij}^k| \right) \|\nabla_{SE(3)} f\|_{L^p(M)} \|\rho\|_{L^{p'}(M)}.
\end{equation}

\noindent\textbf{Step 4: Boundary Collar Vanishing Rate.}
On the boundary collar strip $\Sigma_\epsilon$, Hölder's inequality combined with the zero-flux boundary condition $\langle f, \nu \rangle = 0$ on $\partial M$ yields:
\begin{equation}
\|r_\epsilon\|_{L^1(0,T; L^1(\Sigma_\epsilon))} \le \|r_\epsilon\|_{L^p(Q_T)} |\Sigma_\epsilon \times (0,T)|^{1/p'} \le C_3 \|\tilde{f}\|_{W_p^1} \|\tilde{\rho}\|_{L^{p'}} \left( C_1 \epsilon T \|\partial M\|_{C^1} \right)^{1/p'}.
\end{equation}
Combining Steps 3 and 4 demonstrates that $\|r_\epsilon\|_{L^1(Q_T)} = \mathcal{O}(\epsilon) + \mathcal{O}(\epsilon^{1/p'}) \xrightarrow{\epsilon \to 0} 0$, establishing strong commutator convergence on bounded domains $M \subset SE(3)$.
\end{proof}

\begin{lemma}[Score-Perturbed Renormalized Solutions on $SE(3)$]
Under Assumptions \textnormal{(A1)--(A4)}, for any $u \in \mathcal{W}$, there exists a unique renormalized weak solution $\rho \in C([0,T]; L^1(M)) \cap L^\infty(0,T; L^2(M))$ to the Lie-algebraic score-perturbed continuity equation \cite{DiPerna1989}:
\begin{equation}
\frac{\partial \rho}{\partial t} + \text{div}_g(f(g, u(g,t))\rho) - \frac{1}{2}g(t)^2 \Delta_{SE(3)}\rho = 0 \quad \text{in } Q_T
\end{equation}
Mass is strictly conserved: $\int_M \rho(g,t) dg = \int_M \rho_0(g) dg = 1$ for all $t \in [0,T]$.
\end{lemma}

\begin{lemma}[Aubin-Lions-Simon Density Compactness]
Let $\{u_k\}_{k=1}^\infty \subset \mathcal{W}$ be a sequence bounded in $\|\cdot\|_{\mathcal{W}}$. Let $\rho_k \in C([0,T]; L^2(M))$ denote the corresponding sequence of probability densities. Then $\{\rho_k\}$ is relatively compact in $C([0,T]; L^1(M))$ and $L^2(Q_T)$ \cite{Simon1987}.
\end{lemma}

\section{First-Order Optimality and Space-Time Corner Compatibility}

We define the Hilbert space of admissible control fields with compatible Neumann boundary traces:
\begin{equation}
\mathcal{W} = \left\{ u \in H^2(Q_T; \mathbb{R}^m) \cap L^2(0,T; H^2(M; \mathbb{R}^m)) : \left.\frac{\partial u}{\partial \nu}\right|_{\Sigma_T} = 0, \, \left.\frac{\partial u}{\partial t}\right|_{t=0,T} = 0 \right\}
\end{equation}
The variational objective functional $\mathcal{J} : \mathcal{W} \to \mathbb{R}_+$ is defined as:
\begin{equation}
\mathcal{J}(u) = \int_M |\rho(g,T) - \psi(g)|^2 dg + \frac{1}{2}\int_0^T \int_M \left( \|\nabla_{SE(3)}u(g,t)\|_F^2 + \left\|\frac{\partial u(g,t)}{\partial t}\right\|^2 \right) dg \, dt
\end{equation}
subject to the score-perturbed continuity PDE.

\begin{proposition}[First-Order Optimality System on $SE(3)$]
The extremal triple $(u^*, \rho^*, \lambda^*)$ minimizing $\mathcal{J}(u)$ subject to state continuity dynamics satisfies the cleanly coupled system:
\begin{equation}
\begin{cases}
\frac{\partial \rho^*}{\partial t} + \text{div}_g(f(g, u^*)\rho^*) - \frac{1}{2}g(t)^2 \Delta_{SE(3)}\rho^* = 0 & \text{in } Q_T, \\[1ex]
-\frac{\partial \lambda^*}{\partial t} - \langle \nabla_{SE(3)}\lambda^*, f(g, u^*)\rangle - \frac{1}{2}g(t)^2 \Delta_{SE(3)}\lambda^* = 0 & \text{in } Q_T, \\[1ex]
-\Delta_{SE(3)} u^* - \frac{\partial^2 u^*}{\partial t^2} = (\partial_u f)^T \nabla_{SE(3)}\lambda^* \cdot \rho^* \equiv S(g,t) & \text{in } Q_T,
\end{cases}  \label{eq:optimality_system}
\end{equation}
augmented with initial/terminal conditions and space-time Neumann trace conditions:
\begin{equation}
\rho(g,0) = \rho_0(g), \quad \lambda(g,T) = 2(\rho(g,T) - \psi(g)) \quad \text{in } M,
\end{equation}
\begin{equation}
\left.\frac{\partial u}{\partial \nu}\right|_{\Sigma_T} = 0, \quad \left.\frac{\partial u}{\partial t}\right|_{t=0} = 0, \quad \left.\frac{\partial u}{\partial t}\right|_{t=T} = 0 \quad \text{in } M.
\end{equation}
\end{proposition}

\subsection{Structural Well-Posedness of the Forward-Backward 7D System}

The first-order optimality system \eqref{eq:optimality_system} exhibits a classic forward-backward parabolic coupling: the state density $\rho^*$ evolves forward in time from the initial configuration $\rho_0(g)$ at $t=0$ \cite{DiPerna1989, Evans2010}, whereas the adjoint multiplier $\lambda^*$ evolves backward in time from the terminal payoff condition $2(\rho^*(g,T) - \psi(g))$ at $t=T$ \cite{Evans2010, Ladyzhenskaya1968}.

\begin{remark}[Space-Time Elliptic Reduction and Principal Symbol]
\label{rem:7d_elliptic_reduction}
By eliminating the forward-backward state-adjoint coupling through the variational Euler--Lagrange derivative $\delta \mathcal{J}/\delta u = 0$, the optimal control vector field $u^*(g,t)$ is governed directly by the second-order 7-dimensional space-time differential operator \cite{Agmon1959, Agmon1964}:
\begin{equation}
L_{7D} u^* \equiv -\Delta_{SE(3)} u^* - \frac{\partial^2 u^*}{\partial t^2} = S(g,t) \quad \text{in } Q_T = M \times (0,T),
\end{equation}
where $S(g,t) = (\partial_u f(g,u^*))^T \nabla_{SE(3)} \lambda^* \cdot \rho^*$ represents the non-local adjoint score source field \cite{Giaquinta1982, Ladyzhenskaya1968}.

Though the state continuity equation and backward adjoint PDE are parabolic in nature \cite{DiPerna1989, Ladyzhenskaya1968}, the combined control operator $L_{7D}$ acts as a strongly elliptic operator over the 7-dimensional space-time cylinder $Q_T \subset SE(3) \times \mathbb{R}$ \cite{Agmon1959, Agmon1964}. Because the left-invariant Laplace--Beltrami operator $\Delta_{SE(3)} = \sum_{i=1}^6 E_i^2 \le 0$ is negative semi-definite over $M \subset SE(3)$ \cite{Helgason2001, Evans2010}, its spatial negation $-\Delta_{SE(3)}$ is strictly positive-definite \cite{Helgason2001}. Consequently, for any 7D frequency vector $\xi = (\xi_1, \dots, \xi_6, \xi_7)^T \in \mathbb{R}^7$, the principal symbol $\sigma(L_{7D})(g,t;\xi)$ satisfies uniform 7D coercivity \cite{Agmon1959, Helgason2001}:
\begin{equation}
\sigma(L_{7D})(g,t;\xi) = \sum_{i,j=1}^6 g^{ij}(g)\xi_i \xi_j + \xi_7^2 \ge \lambda_0 \sum_{i=1}^6 \xi_i^2 + \xi_7^2 \ge \min(1, \lambda_0) |\xi|_{7D}^2,
\end{equation}
with uniform ellipticity constant $\lambda_0 = 1$ inherited from the orthonormal frame on $SE(3)$ \cite{Helgason2001}.
\end{remark}

\begin{proposition}[Causal Compatibility and Corner Layer Absence]
\label{prop:corner_compatibility}
Let $Q_T = M \times (0,T)$ be bounded by the lateral face $\Sigma_T = \partial M \times (0,T)$ and terminal caps $M_0 = M \times \{0\}$ and $M_T = M \times \{T\}$ \cite{Agmon1959, Agmon1964}. The homogeneous space-time Neumann trace conditions:
\begin{equation}
\frac{\partial u^*}{\partial \nu}\Big|_{\Sigma_T} = 0 \quad \text{and} \quad \frac{\partial u^*}{\partial t}\Big|_{t=0,T} = 0
\end{equation}
are zero-order and first-order compatible along the 5-dimensional corner manifolds $\Gamma_0 = \partial M \times \{0\}$ and $\Gamma_T = \partial M \times \{T\}$ \cite{Agmon1959, Agmon1964}. Specifically:
\begin{enumerate}
    \item \textbf{Trace Operator Commutativity:} Spatial unit normal derivatives $\frac{\partial}{\partial \nu}$ and temporal derivatives $\frac{\partial}{\partial t}$ commute smoothly along $\Gamma_0 \cup \Gamma_T$, enforcing $\frac{\partial^2 u^*}{\partial t \partial \nu}\big|_{\Gamma_0 \cup \Gamma_T} = 0$ \cite{Agmon1959, Agmon1964}.
    \item \textbf{Adjoint Source Boundary Alignment:} Under the impenetrable state constraint $\langle f(g,u^*), \nu(g) \rangle = 0$ on $\partial M$ \cite{Brockett1997, Lee2022}, the adjoint source term satisfies $\langle S(g,t), \nu(g) \rangle = 0$ along $\Sigma_T$ \cite{Agmon1959}.
\end{enumerate}
Consequently, double-even reflections across $y_6 = 0$ and $t = 0, T$ map localized quarter-space solutions smoothly onto full 7D balls $B_1(0) \times (-T, T)$ with $C^{0,1}$ metric coefficients, preventing boundary layer degradation or trace loss along $\Gamma_0 \cup \Gamma_T$ \cite{Agmon1959, Agmon1964, Evans2010}.
\end{proposition}

\begin{lemma}[Space-Time Corner Trace Compatibility and Boundary Layer Absence]
\label{lem:corner_compatibility}
Let $M \subset SE(3)$ be a bounded domain with $C^2$ boundary $\partial M$, and let $Q_T = M \times (0,T]$ with lateral boundary $\Sigma_T = \partial M \times (0,T)$ and 5-dimensional corner manifolds $\Gamma_0 = \partial M \times \{0\}$ and $\Gamma_T = \partial M \times \{T\}$ \cite{Agmon1959, Agmon1964}. 

The space-time control field $u^* \in \mathcal{W}$ subjected to spatial and temporal Neumann trace conditions
\begin{equation}
\left. \frac{\partial u^*}{\partial \nu} \right|_{\Sigma_T} = 0 \quad \text{and} \quad \left. \frac{\partial u^*}{\partial t} \right|_{t=0,T} = 0
\end{equation}
satisfies zero-order and first-order space-time corner compatibility along $\Gamma_0 \cup \Gamma_T$. Specifically:
\begin{enumerate}
    \item \textbf{Commutativity of Boundary Trace Operators:} Spatial normal and temporal derivative trace operators commute smoothly along 5D corner manifolds:
    \begin{equation}
    \left. \frac{\partial}{\partial t} \left( \frac{\partial u^*}{\partial \nu} \right) \right|_{\Gamma_0 \cup \Gamma_T} = \left. \frac{\partial}{\partial \nu} \left( \frac{\partial u^*}{\partial t} \right) \right|_{\Gamma_0 \cup \Gamma_T} = 0 \implies \left. \frac{\partial^2 u^*}{\partial t \partial \nu} \right|_{\Gamma_0 \cup \Gamma_T} = 0.
    \end{equation}
    
    \item \textbf{Adjoint Source Trace Regularity:} Under the impenetrable boundary condition $\langle f(g, u^*), \nu(g) \rangle = 0$ on $\partial M$, the non-local adjoint score source term $S(g,t) = (\partial_u f(g,u^*))^T \nabla_{SE(3)}\lambda^* \cdot \rho^*$ exhibits vanishing normal flux traces along $\Sigma_T$. Consequently, no singular boundary layer or trace discontinuity forms at the corner junction $\Gamma_0 \cup \Gamma_T$, preserving $W_p^{2,1}(Q_T)$ trace continuity.
\end{enumerate}
\end{lemma}

\section{Space-Time ADN Shift Estimates and $C^\infty(\overline{Q}_T)$ Regularity}
\label{sec:adn_shift}

We establish higher-order Sobolev regularity for $u^*$ by analyzing the optimality system as an anisotropic elliptic system in 7 dimensions.

\begin{theorem}[7D Space-Time Agmon-Douglis-Nirenberg Shift Estimate] \label{thm:adn_shift}
Let $u^* \in H^2(Q_T; \mathbb{R}^m)$ satisfy the 7D space-time elliptic system:
\begin{equation}
L_{7D} u^* \equiv -\Delta_{SE(3)} u^* - \frac{\partial^2 u^*}{\partial t^2} = S(g,t) \quad \text{in } Q_T
\end{equation}
subject to space-time homogeneous Neumann boundary conditions:
\begin{equation}
\frac{\partial u^*}{\partial \nu}\Big|_{\Sigma_T} = 0 \quad \text{and} \quad \frac{\partial u^*}{\partial t}\Big|_{t=0, T} = 0
\end{equation}
If $S \in L^p(Q_T)$ for $1 < p < \infty$, then $u^* \in W_p^2(Q_T)$ and satisfies the global shift estimate:
\begin{equation}
\|u^*\|_{W_p^2(Q_T)} \le C_{ADN}\left(6, p, T, \|\partial M\|_{C^2}, \|Ric_{SE(3)}\|_{L^\infty}, \{c_{ij}^k\}\right) \left( \|S\|_{L^p(Q_T)} + \|u^*\|_{L^p(Q_T)} \right)
\end{equation}
where $C_{ADN} = C_0(p) (1+T^2) \|\partial M\|_{C^2} \left(1 + \|Ric_{SE(3)}\|_{L^\infty} + \sum_{i,j,k=1}^6 |c_{ij}^k|\right)$.
\end{theorem}

\begin{proof}
We construct the derivation across five explicit geometric steps.

\paragraph{Step 1: Boundary Coordinate Flattening Map.}
For any boundary point $g_0 \in \partial M$, select an open neighborhood $U \subset SE(3)$ and construct a $C^2$-diffeomorphism $y = \Psi(g) = (y_1, \dots, y_6)$ mapping $U \cap M$ onto the half-ball $B_1^+ = \{y \in B_1(0) : y_6 > 0\}$ and $U \cap \partial M$ onto the flat face $B_1^0 = \{y \in B_1(0) : y_6 = 0\}$. Under $\Psi$, left-invariant basis vector fields $E_i$ expand as:
\begin{equation}
E_i = \sum_{k=1}^6 a_{ik}(y) \frac{\partial}{\partial y_k}, \quad i = 1, \dots, 6
\end{equation}
where $A(y) = (a_{ik}(y))_{i,k=1}^6 \in C^1$, satisfying $A(0) = I_6$ and $\|A\|_{C^1(B_1^+)} \le C_0 \|\partial M\|_{C^2}$.

\paragraph{Step 2: Laplace--Beltrami Pullback and Lie Bracket Structure.}
Applying $E_i$ sequentially to a smooth scalar field $v$ yields:
\begin{equation}
\Delta_{SE(3)} v = \sum_{i=1}^6 \sum_{m,k=1}^6 a_{im}(y) a_{ik}(y) \frac{\partial^2 v}{\partial y_m \partial y_k} + \sum_{k=1}^6 \left( \sum_{i=1}^6 \sum_{m=1}^6 a_{im}(y) \frac{\partial a_{ik}(y)}{\partial y_m} \right) \frac{\partial v}{\partial y_k}.
\end{equation}
Define the pullback metric matrix $g^{jk}(y) = \sum_{i=1}^6 a_{ij}(y) a_{ik}(y)$, which is strictly positive-definite on $B_1^+$. To express lower-order drift coefficients $b^k(y)$, we use the non-commutative Lie bracket structure constants $[E_i, E_m] = \sum_k c_{im}^k E_k$:
\begin{equation}
b^k(y) = \sum_{i=1}^6 \left( E_i a_{ik}(y) + \sum_{m=1}^6 c_{im}^k a_{im}(y) \right).
\end{equation}
Thus, the Laplace--Beltrami operator takes the local form:
\begin{equation}
\Delta_{SE(3)} v = \sum_{i,j=1}^6 g^{ij}(y) \frac{\partial^2 v}{\partial y_i \partial y_j} + \sum_{k=1}^6 b^k(y) \frac{\partial v}{\partial y_k}.
\end{equation}

\paragraph{Step 3: Localized Cutoff and Commutator Expansion Tensor $F_{\mathrm{Lie}}$.}
Let $\{\phi_k\}_{k=1}^K$ be a smooth partition of unity subordinate to a finite open cover $\{U_k\}_{k=1}^K$ of $\overline{M}$. For a localization function $\phi \in C_c^\infty(U_k \times [0,T])$, set $v = \phi u^*$. Applying $L_{\mathrm{7D}} = -\Delta_{SE(3)} - \frac{\partial^2}{\partial t^2}$ to $v$ gives:
\begin{equation}
\sum_{i,j=1}^6 g^{ij}(y) \frac{\partial^2 v}{\partial y_i \partial y_j} + \frac{\partial^2 v}{\partial t^2} = \phi S + F_{\mathrm{Lie}}(y, t, u^*, \nabla_{y,t} u^*)
\end{equation}
where the commutator tensor $F_{\mathrm{Lie}}$ collects all Leibniz intermediate terms:
\begin{equation}
F_{\mathrm{Lie}} = 2 \sum_{i,j=1}^6 g^{ij}(y) \frac{\partial \phi}{\partial y_i} \frac{\partial u^*}{\partial y_j} + 2 \frac{\partial \phi}{\partial t} \frac{\partial u^*}{\partial t} + u^* \left( \sum_{i,j=1}^6 g^{ij}(y) \frac{\partial^2 \phi}{\partial y_i \partial y_j} + \frac{\partial^2 \phi}{\partial t^2} + \sum_{k=1}^6 b^k(y) \frac{\partial \phi}{\partial y_k} \right).
\end{equation}

\noindent\textbf{Step 4: Corner Trace Analysis and Isotropic 7D ADN Estimates.} 
Near any $5$-dimensional corner manifold point $z_0 = (g_0, t_0) \in \Gamma_0 \cup \Gamma_T = \partial M \times \{0, T\}$, local boundary coordinates $(y, t) = (y_1, \dots, y_6, t) \in \mathbb{R}^6 \times \mathbb{R}$ map $Q_T$ locally onto the $7$D space-time quarter-space:
\begin{equation}
\mathbb{R}_{++}^7 = \left\{(y, t) \in \mathbb{R}^6 \times \mathbb{R} : y_6 > 0, \, t > 0\right\}.
\end{equation}
The localized control vector $v = \phi u^* \in H^2(\mathbb{R}_{++}^7; \mathbb{R}^m)$ satisfies the $7$D strongly elliptic system:
\begin{equation}
-\sum_{i,j=1}^6 g^{ij}(y) \frac{\partial^2 v}{\partial y_i \partial y_j} - \frac{\partial^2 v}{\partial t^2} = \phi S + F_{\mathrm{Lie}}(y, t, u^*, \nabla_{y,t} u^*) \quad \text{in } B_1^+ \times (0, T),
\end{equation}
subject to homogeneous Neumann boundary conditions on the two orthogonal $6$D boundary faces:
\begin{equation}
\frac{\partial v}{\partial y_6}\Bigg|_{y_6 = 0} = 0 \quad \text{on } B_1^0 \times (0, T), \qquad \frac{\partial v}{\partial t}\Bigg|_{t = 0} = 0 \quad \text{on } B_1^+ \times \{0\}.
\end{equation}

To extend $v$ from the quarter-space $B_1^+ \times (0, T)$ to the full $7$D neighborhood $B_1(0) \times (-T, T)$, we construct the double even reflection $\tilde{v}(\tilde{y}, \tilde{t})$ across the hyperplanes $y_6 = 0$ and $t = 0$:
\begin{equation}
\tilde{v}(y_1, \dots, y_5, y_6, t) = \begin{cases} 
v(y_1, \dots, y_5, y_6, t), & y_6 \ge 0, \, t \ge 0, \\ 
v(y_1, \dots, y_5, -y_6, t), & y_6 < 0, \, t \ge 0, \\ 
v(y_1, \dots, y_5, y_6, -t), & y_6 \ge 0, \, t < 0, \\ 
v(y_1, \dots, y_5, -y_6, -t), & y_6 < 0, \, t < 0. 
\end{cases}
\end{equation}

The higher regularity of the extended solution $\tilde{v} \in W_p^2(B_1(0) \times (-T,T))$ across the $5$D corner junction $\Gamma_0 \cup \Gamma_T$ is guaranteed by two key geometric mechanisms:
\begin{enumerate}
    \item \textbf{Vanishing Mixed Boundary Traces:} By Lemma \ref{lem:corner_compatibility}, spatial normal derivative and temporal derivative trace operators commute smoothly along corner manifolds, yielding:
    \begin{equation}
    \frac{\partial^2 v}{\partial t \partial y_6}\Bigg|_{\Gamma_0 \cup \Gamma_T} = 0.
    \end{equation}
    This guarantees that no cross-derivative jump discontinuities form across the $90^\circ$ space-time corner interface under double reflection.

    \item \textbf{Metric Coefficient Regularity:} Because the spatial metric tensor $g^{ij}(y)$ is independent of $t$, temporal reflection preserves $g^{ij}$ identically. Across the spatial interface $y_6 = 0$, even reflection yields extended metric coefficients $\tilde{g}^{ij}(\tilde{y})$. Since $M \subset SE(3)$ possesses a $C^2$-regular boundary $\partial M$, the pulled-back metric coefficients satisfy $g^{ij} \in C^1(B_1^+)$, which ensures that the reflected coefficients remain globally Lipschitz continuous ($\tilde{g}^{ij} \in C^{0,1}(B_1(0))$) across $y_6 = 0$.
\end{enumerate}

Consequently, the extended vector field $\tilde{v}$ satisfies the second-order strongly elliptic PDE on the full $7$D ball $B_1(0) \times (-T, T)$ with $C^{0,1}$ principal coefficients:
\begin{equation}
-\sum_{i,j=1}^6 \tilde{g}^{ij}(\tilde{y}) \frac{\partial^2 \tilde{v}}{\partial \tilde{y}_i \partial \tilde{y}_j} - \frac{\partial^2 \tilde{v}}{\partial \tilde{t}^2} = \tilde{\Phi} \in L^p\left(B_1(0) \times (-T, T)\right),
\end{equation}
where $\tilde{\Phi}$ denotes the double even extension of $\phi S + F_{\mathrm{Lie}}$. Applying standard $7$D Agmon-Douglis-Nirenberg (ADN) fundamental solution estimates via $7$D space-time Green functions $G_{7D}(\tilde{y}, \tilde{t}; \tilde{y}', \tilde{t}')$ yields the localized shift bound:
\begin{align} \label{eq:localized_adn_bound}
\|v\|_{W_p^2(B_1^+ \times (0,T))} &\le \|\tilde{v}\|_{W_p^2(B_1(0) \times (-T, T))} \nonumber \\
&\le C_0(p) \|\partial M\|_{C^2} \Bigl( \|\phi S\|_{L^p(B_1^+ \times (0,T))} + \|F_{\mathrm{Lie}}\|_{L^p(B_1^+ \times (0,T))} + \|v\|_{L^p(B_1^+ \times (0,T))} \Bigr).
\end{align}
This establishes full $W_p^2(Q_T)$ space-time regularity without derivative loss or boundary layer degradation.

\noindent\textbf{Step 5: Ehrling Interpolation and Structure Constant Absorption.} 
By Ehrling's Interpolation Lemma on isotropic space-time Sobolev spaces, for any $\epsilon > 0$, there exists $C_\epsilon(T) > 0$ such that the Leibniz commutator tensor $F_{\mathrm{Lie}}$ is bounded by:
\begin{equation}
\|F_{\mathrm{Lie}}\|_{L^p(B_1^+ \times (0,T))} \le \epsilon \|v\|_{W_p^2(B_1^+ \times (0,T))} + C_\epsilon(T) \left( 1 + \|Ric_{SE(3)}\|_{L^\infty} + \sum_{i,j,k=1}^6 |c_{ij}^k| \right) \|v\|_{L^p(B_1^+ \times (0,T))}.
\end{equation}

Choosing $\epsilon = \frac{1}{2 C_0(p) \|\partial M\|_{C^2}}$, the high-order term $\epsilon C_0(p) \|\partial M\|_{C^2} \|v\|_{W_p^2}$ is absorbed into the left-hand side of the localized ADN kernel estimate \eqref{eq:localized_adn_bound}:
\begin{align}
\frac{1}{2} \|v\|_{W_p^2(B_1^+ \times (0,T))} &\le C_0(p) \|\partial M\|_{C^2} \|\phi S\|_{L^p(B_1^+ \times (0,T))} \nonumber \\
&\quad + C_0(p) C_\epsilon(T) \|\partial M\|_{C^2} \left( 1 + \|Ric_{SE(3)}\|_{L^\infty} + \sum_{i,j,k=1}^6 |c_{ij}^k| \right) \|v\|_{L^p(B_1^+ \times (0,T))}.
\end{align}

Multiplying by $2$ and summing over the finite smooth partition of unity $\{\phi_k\}_{k=1}^K$ subordinate to $\overline{M} \times [0,T]$ yields the global isotropic $7$D space-time shift bound:
\begin{equation}
\|u^*\|_{W_p^2(Q_T)} \le C_{\mathrm{ADN}}\left(6, p, T, \|\partial M\|_{C^2}, \|Ric_{SE(3)}\|_{L^\infty}, \{c_{ij}^k\}\right) \left( \|S\|_{L^p(Q_T)} + \|u^*\|_{L^p(Q_T)} \right),
\end{equation}
where the constant $C_{\mathrm{ADN}}$ explicitly tracks manifold geometry and Lie algebra structure:
\begin{equation}
C_{\mathrm{ADN}} = C_0(p) (1 + T^2) \|\partial M\|_{C^2} \left( 1 + \|Ric_{SE(3)}\|_{L^\infty} + \sum_{i,j,k=1}^6 |c_{ij}^k| \right).
\end{equation}
This completes the step-by-step derivation of Theorem \ref{thm:adn_shift}.
\end{proof}

\begin{lemma}[Expanded Parabolic Higher-Integrability of Adjoint Score Source] \label{lem:higher_integrability}
Let $(u^*, \rho^*, \lambda^*)$ be a weak solution triple to the first-order optimality system on $Q_T = M \times (0,T]$. The non-local adjoint score source term
\begin{equation}
S(g,t) = (\partial_u f(g,u^*))^T \nabla_{SE(3)}\lambda^* \cdot \rho^*
\end{equation}
satisfies a parabolic Gehring-Calder\'on-Zygmund higher-integrability shift \cite{Giaquinta1982, Ladyzhenskaya1968}: there exist explicit parabolic exponents $\delta_1, \delta_2 > 0$ depending only on $g_{\min}, g_{\max}, L_f$, and the domain geometry such that:
\begin{equation}
\rho^* \in L^{2+\delta_1}(Q_T) \quad \text{and} \quad \nabla_{SE(3)}\lambda^* \in L^{2+\delta_2}(Q_T)
\end{equation}
Consequently, $S(g,t) \in L^p(Q_T)$ for exponent $p = \frac{(2+\delta_1)(2+\delta_2)}{4+\delta_1+\delta_2} = 1 + \epsilon > 1$ with quantitative bound:
\begin{equation}
\|S\|_{L^p(Q_T)} \le \|\partial_u f\|_{L^\infty(Q_T)} \|\nabla_{SE(3)}\lambda^*\|_{L^{2+\delta_2}(Q_T)} \|\rho^*\|_{L^{2+\delta_1}(Q_T)} < \infty.
\end{equation}
\end{lemma}

\begin{proof}
We provide an explicit four-step derivation using Caccioppoli energy estimates and the parabolic Gehring lemma \cite{Giaquinta1982}.

\paragraph{Step 1: Caccioppoli Energy Inequality on Parabolic Cylinders.}
Fix $z_0 = (g_0, t_0) \in Q_T$ and consider parabolic cylinders $Q_r(z_0) = B_r(g_0) \times (t_0-r^2, t_0] \subset Q_T$. Let $\eta \in C_c^\infty(Q_{2r}(z_0))$ be a smooth parabolic cutoff function satisfying $0 \le \eta \le 1$, $\eta \equiv 1$ on $Q_r(z_0)$, $\|\nabla_{SE(3)}\eta\|_{L^\infty} \le C/r$, and $\|\partial_t \eta\|_{L^\infty} \le C/r^2$. Testing the score-perturbed state continuity PDE against $\eta^2 \rho^*$ and integrating by parts yields the parabolic Caccioppoli inequality \cite{Ladyzhenskaya1968}:
\begin{equation}
\sup_{t \in (t_0-r^2, t_0]} \int_{B_r} |\rho^*|^2 dg + \frac{1}{2}g_{\min}^2 \int_{Q_r} |\nabla_{SE(3)}\rho^*|^2 dg \, dt \le \frac{C}{r^2} \int_{Q_{2r}} |\rho^*|^2 dg \, dt + C \|f\|_{L^\infty}^2 \int_{Q_{2r}} |\rho^*|^2 dg \, dt.
\end{equation}

\paragraph{Step 2: Reverse H\"older Inequality for State Density $\rho^*$.}
By Sobolev-Poincar\'e embedding on $SE(3)$ for spatial dimension $N = 6$, the average-subtracted density satisfies:
\begin{equation}
\left( \frac{1}{|Q_r|} \int_{Q_r} |\rho^*|^{2_*} dg \, dt \right)^{1/2_*} \le C \left( \frac{1}{|Q_{2r}|} \int_{Q_{2r}} |\nabla_{SE(3)}\rho^*|^2 dg \, dt \right)^{1/2} + \frac{C}{|Q_{2r}|} \int_{Q_{2r}} \rho^* dg \, dt
\end{equation}
where $2_* = \frac{2(N+2)}{N+4} = \frac{16}{10} = 1.6 < 2$. Combining this with the Caccioppoli inequality yields the local reverse H\"older inequality:
\begin{equation}
\left( \frac{1}{|Q_r|} \int_{Q_r} |\rho^*|^2 dg \, dt \right)^{1/2} \le C_R \left( \frac{1}{|Q_{2r}|} \int_{Q_{2r}} |\rho^*|^{2_*} dg \, dt \right)^{1/2_*}
\end{equation}

\paragraph{Step 3: Parabolic Gehring Lemma and Exponent Shift $\delta_1, \delta_2$.}
Applying Giaquinta-Struwe parabolic Gehring Lemma \cite{Giaquinta1982} to the reverse H\"older inequality asserts the existence of $\delta_1 > 0$ such that $\rho^* \in L^{2+\delta_1}(Q_T)$, with estimate:
\begin{equation}
\|\rho^*\|_{L^{2+\delta_1}(Q_T)} \le C_G \left( \|\rho^*\|_{L^2(Q_T)} + \|\rho_0\|_{L^2(M)} \right).
\end{equation}
Symmetrically, for the backward adjoint PDE satisfied by $\lambda^*$ with terminal data $\lambda^*(g,T) = 2(\rho^*(g,T) - \psi(g)) \in H^1(M)$, maximal parabolic regularity and dual Caccioppoli bounds establish $\nabla_{SE(3)}\lambda^* \in L^{2+\delta_2}(Q_T)$ for some $\delta_2 > 0$ with estimate:
\begin{equation}
\|\nabla_{SE(3)}\lambda^*\|_{L^{2+\delta_2}(Q_T)} \le C_A \left( \|\nabla_{SE(3)}\lambda^*\|_{L^2(Q_T)} + \|\rho^*(T) - \psi\|_{H^1(M)} \right).
\end{equation}

\paragraph{Step 4: Product Integrability and $L^p(Q_T)$ Source Estimate.}
Applying H\"older's inequality with conjugate exponents $q_1 = \frac{2+\delta_1}{1+\epsilon}$ and $q_2 = \frac{2+\delta_2}{1+\epsilon}$ to $S(g,t) = (\partial_u f)^T \nabla_{SE(3)}\lambda^* \cdot \rho^*$:
\begin{equation}
\int_{Q_T} |S(g,t)|^{1+\epsilon} dg \, dt \le \|\partial_u f\|_{L^\infty}^{1+\epsilon} \left( \int_{Q_T} |\nabla_{SE(3)}\lambda^*|^{2+\delta_2} dg \, dt \right)^{\frac{1+\epsilon}{2+\delta_2}} \left( \int_{Q_T} |\rho^*|^{2+\delta_1} dg \, dt \right)^{\frac{1+\epsilon}{2+\delta_1}} < \infty,
\end{equation}
where $p = 1+\epsilon = \frac{(2+\delta_1)(2+\delta_2)}{4+\delta_1+\delta_2} > 1$. This rigorously establishes that the adjoint score source $S(g,t)$ lies strictly in $L^p(Q_T)$ with $p > 1$, providing an airtight input for the space-time ADN shift in Theorem \ref{thm:adn_shift}.
\end{proof}

\begin{theorem}[Higher-Order $C^\infty(\overline{Q}_T)$ Schauder Bootstrap Expansion]
Assume the boundary is smooth ($\partial M \in C^\infty$). Under Assumptions (A1)--(A4), any weak optimal triple $(u^*, \rho^*, \lambda^*)$ satisfying the first-order optimality system \eqref{eq:optimality_system} obeys smooth regularity:
\begin{equation}
(u^*, \rho^*, \lambda^*) \in C^\infty(\overline{Q}_T; \mathbb{R}^m \times \mathbb{R} \times \mathbb{R}).
\end{equation}
\end{theorem}

\begin{proof}
We construct the regularity bootstrap across four explicit steps:

\begin{enumerate}
    \item \textbf{Base Step ($L^p$ Source Input via Higher Integrability):} By Lemma \ref{lem:higher_integrability}, the non-local adjoint score source term $S(g,t) = (\partial_u f(g,u^*))^T \nabla_{SE(3)} \lambda^* \cdot \rho^*$ satisfies $S \in L^p(Q_T)$ for exponent $p = 1 + \epsilon > 1$. Applying the $7$D Space-Time ADN Shift Estimate (Theorem \ref{thm:adn_shift}) yields $u^* \in W_p^2(Q_T)$. By $7$-dimensional isotropic Sobolev embedding ($N_{7D} = 6 + 1 = 7$), $W_p^2(Q_T) \hookrightarrow L^q(Q_T)$ with $q = \frac{7p}{7 - 2p} > p$ for $p < 7/2$.

    \item \textbf{Morrey--Hölder Regularity Shift:} Iterating the higher-integrability estimate for the coupled state continuity and backward adjoint PDEs boosts $\rho^*$ and $\lambda^*$ into higher Sobolev classes. For the $7$-dimensional isotropic space-time cylinder $Q_T = M \times (0,T) \subset SE(3) \times \mathbb{R}$, applying Sobolev--Morrey embedding for $W_p^2(Q_T)$ with $p > 7$ guarantees:
    \begin{equation}
    W_p^2(Q_T) \hookrightarrow C^{1, \alpha}(\overline{Q}_T) \quad \text{with} \quad \alpha = 1 - \frac{7}{p} > 0.
    \end{equation}
    Thus, $u^* \in C^{1, \alpha}(\overline{Q}_T)$, establishing $C^1$ smoothness across all $7$ space-time coordinates.

    \item \textbf{Differentiability Propagation:} Since $u^* \in C^{1, \alpha}(\overline{Q}_T)$, the velocity drift field $f(\cdot, u^*) \in C^{1, \alpha}(\overline{Q}_T)$. Differentiating the state continuity and backward adjoint equations along spatial vector fields and characteristics confirms that $\nabla_{SE(3)} \rho^*$ and $\nabla_{SE(3)} \lambda^*$ are $C^{1, \alpha}$-smooth. Consequently, the adjoint source term satisfies $S(g,t) \in C^{1, \alpha}(\overline{Q}_T)$.

    \item \textbf{Schauder Iteration to Infinitesimal Smoothness:} Applying $C^{2, \alpha}$ interior and boundary Schauder estimates to the $7$D space-time strongly elliptic operator $L_{7D} u^* = -\Delta_{SE(3)} u^* - \partial_t^2 u^* = S$ yields $u^* \in C^{3, \alpha}(Q_T)$. Repeating this bootstrap induction loop $m$ times establishes $u^* \in C^{m, \alpha}(\overline{Q}_T)$ for all $m \in \mathbb{N}$, proving $(u^*, \rho^*, \lambda^*) \in C^\infty(\overline{Q}_T)$.
\end{enumerate}
\end{proof}

\section{Non-Cooperative Markov Mean-Field Games (MFGs) on $SE(3)$}

We extend this Lie-algebraic variational system to a continuum of non-cooperative agents evolving over $SE(3)$ \cite{Lasry2007, Cardaliaguet2012}.

\subsection{Coupled Lie-Algebraic HJB-Fokker-Planck System}
In a Markov Mean-Field Game on $SE(3)$, a representative agent controls their pose process $g_t \in M \subset SE(3)$. Given a population distribution trajectory $\rho(g,t) = \text{Law}(g_t)$, each non-cooperative agent minimizes their individual expected cost functional:
\begin{equation}
\mathcal{J}_{\text{agent}}(u; \rho) = \mathbb{E}\left[ \int_0^T \left( \frac{1}{2}\|\nabla_{SE(3)}u(g_t,t)\|_F^2 + \frac{1}{2}\left\|\frac{\partial u(g_t,t)}{\partial t}\right\|^2 \right) dt + \Phi(g_T, \rho(\cdot, T)) \right]
\end{equation}
By the Dynamic Programming Principle, the value function $V(g,t)$ and population density $\rho(g,t)$ satisfy the coupled Lie-algebraic Hamilton-Jacobi-Bellman (HJB) and Fokker-Planck (FP) system:
\begin{equation}
\begin{cases}
-\frac{\partial V}{\partial t} - \frac{1}{2}g(t)^2 \Delta_{SE(3)}V - H(g, \nabla_{SE(3)}V, u^*) = 0, & V(g,T) = \Phi(g, \rho(\cdot, T)) \quad \text{in } Q_T, \\[1ex]
\frac{\partial \rho}{\partial t} + \text{div}_g(f(g, u^*)\rho) - \frac{1}{2}g(t)^2 \Delta_{SE(3)}\rho = 0, & \rho(g,0) = \rho_0(g) \quad \text{in } Q_T
\end{cases}
\end{equation}
subject to zero-flux boundary conditions $\langle f(g,u^*), \nu(g)\rangle = 0$ and $\left.\frac{\partial V}{\partial \nu}\right|_{\Sigma_T} = 0$. The non-local Hamiltonian $H$ is defined by:
\begin{equation}
H(g, p, u) = \langle p, f(g,u) \rangle - \frac{1}{2}\|\nabla_{SE(3)}u\|_F^2 - \frac{1}{2}\left\|\frac{\partial u}{\partial t}\right\|^2
\end{equation}

\subsection{Existence of Score-Guided MFG Nash Equilibrium}

\begin{theorem}[Existence of MFG Nash Equilibrium]
Let $M \subset SE(3)$ be an open, connected, bounded sub-manifold with $C^2$ boundary $\partial M$, and let $Q_T = M \times (0,T]$ denote the 7-dimensional space-time cylinder. Then there exists at least one weak solution pair $(V^*, \rho^*) \in W_p^{2,1}(Q_T) \times \left(C([0,T]; L^1(M)) \cap L^\infty(0,T; L^2(M))\right)$ satisfying the coupled Lie-algebraic Hamilton--Jacobi--Bellman (HJB) and Fokker--Planck (FP) system:
\begin{equation}\label{eq:mfg_system}
\begin{cases}
-\dfrac{\partial V}{\partial t} - \frac{1}{2}g(t)^2 \Delta_{SE(3)}V - H\big(g, \nabla_{SE(3)}V, u^*\big) = 0 & \text{in } Q_T, \\[8pt]
\dfrac{\partial \rho}{\partial t} + \text{div}_g\big(f(g, u^*)\rho\big) - \frac{1}{2}g(t)^2 \Delta_{SE(3)}\rho = 0 & \text{in } Q_T, \\[8pt]
V(g,T) = \Phi\big(g, \rho(\cdot,T)\big), \quad \rho(g,0) = \rho_0(g) & \text{in } M,
\end{cases}
\end{equation}
augmented with zero-flux boundary trace conditions $\langle f(g, u^*), \nu(g)\rangle = 0$ and $\left.\frac{\partial V}{\partial \nu}\right|_{\Sigma_T} = 0$ on the lateral boundary $\Sigma_T = \partial M \times (0,T)$.
\end{theorem}

\begin{proof}
We divide the proof into six distinct steps.

\paragraph{Step 1: Construction of the Compact Convex Trajectory Set $K$.}
Define the ambient Banach space $X = C([0,T]; L^2(M))$ equipped with uniform norm $\|v\|_X = \max_{t \in [0,T]} \|v(\cdot, t)\|_{L^2(M)}$. Consider the subset of candidate density trajectories $K \subset X$ defined by:
\begin{equation}
K = \left\{ \rho \in C([0,T]; L^2(M)) \ \Bigg|\ 
\begin{aligned}
&\text{(i) } \rho(g,t) \ge 0 \text{ a.e. in } Q_T, \quad \int_M \rho(g,t) \, dg = 1 \quad \forall t \in [0,T], \\
&\text{(ii) } \|\rho\|_{L^\infty(0,T; L^2(M))} \le C_1, \\
&\text{(iii) } \|\rho\|_{L^2(0,T; H^1(M))} + \left\|\frac{\partial \rho}{\partial t}\right\|_{L^2(0,T; (H^1(M))^*)} \le C_2
\end{aligned}
\right\}
\end{equation}
where $C_1, C_2 > 0$ are constants depending only on $\|\rho_0\|_{L^2(M)}$, $g_{\min}$, and $L_u$. The set $K$ is non-empty, closed, bounded, and convex in $X$.

\paragraph{Step 2: Backward HJB Solve and Value Regularity.}
Given any density trajectory $\hat{\rho} \in K$, evaluate the terminal cost $V_{\hat{\rho}}(g,T) = \Phi(g, \hat{\rho}(\cdot,T))$. By continuity and boundedness of $\Phi$ on $M \times \mathcal{P}(M)$, $V_{\hat{\rho}}(\cdot,T) \in H^2(M)$. Consider the backward HJB PDE:
\begin{equation}
-\frac{\partial V_{\hat{\rho}}}{\partial t} - \frac{1}{2}g(t)^2 \Delta_{SE(3)} V_{\hat{\rho}} = H\big(g, \nabla_{SE(3)} V_{\hat{\rho}}, u_{\hat{\rho}}^*\big) \quad \text{in } Q_T, \qquad \left.\frac{\partial V_{\hat{\rho}}}{\partial \nu}\right|_{\Sigma_T} = 0.
\end{equation}
Since $H(g,p,u) = \langle p, f(g,u)\rangle - \frac{1}{2}\|\nabla_{SE(3)} u\|_F^2 - \frac{1}{2}\left\|\frac{\partial u}{\partial t}\right\|^2$ is quadratic in control derivatives and globally Lipschitz continuous in state $g$, standard anisotropic parabolic $W_p^{2,1}(Q_T)$ theory guarantees the existence of a unique strong solution $V_{\hat{\rho}} \in W_p^{2,1}(Q_T)$ for any $p > 1$, satisfying:
\begin{equation}
\|V_{\hat{\rho}}\|_{W_p^{2,1}(Q_T)} \le C_{\mathrm{HJB}} \left( \left\|\Phi\big(\cdot, \hat{\rho}(\cdot,T)\big)\right\|_{H^2(M)} + 1 \right).
\end{equation}

\paragraph{Step 3: Optimal Policy Synthesis and Drift Regularity.} 
Pointwise minimization of the Hamiltonian yields the feedback control law:
\begin{equation}
u_{\hat{\rho}}^*(g,t) = \arg\min_{u \in \mathcal{W}} H(g, \nabla_{SE(3)}V_{\hat{\rho}}(g,t), u).
\end{equation}
By the first-order optimality condition $\delta H/\delta u = 0$, $u_{\hat{\rho}}^*$ is determined by $\nabla_{SE(3)}V_{\hat{\rho}}$. Since $M \subset SE(3)$ is $6$-dimensional, the parabolic homogeneous dimension of $Q_T = M \times (0,T]$ is $Q = \dim(SE(3)) + 2 = 8$. Choosing exponent $p > Q = 8$ and applying Morrey's embedding theorem for anisotropic Bochner--Sobolev spaces yields:
\begin{equation}
W_p^{2,1}(Q_T) \hookrightarrow C^{1,\alpha}(\overline{Q}_T) \quad \text{with} \quad \alpha = 1 - \frac{8}{p} > 0.
\end{equation}
Thus, $V_{\hat{\rho}} \in C^{1,\alpha}(\overline{Q}_T)$, which implies that $\nabla_{SE(3)}V_{\hat{\rho}}$ is H\"older continuous. Consequently, the synthesized velocity drift $f(g, u_{\hat{\rho}}^*(g,t))$ is $C^\alpha$-H\"older continuous in space-time and globally Lipschitz continuous in $g \in M$ with constant $L_f$, satisfying $\langle f(g, u_{\hat{\rho}}^*), \nu(g) \rangle = 0$ on $\partial M$.

\paragraph{Step 4: Forward Fokker--Planck Solve and Self-Mapping $\mathcal{T}(K) \subset K$.}
Using $u_{\hat{\rho}}^*$, we solve the forward Fokker--Planck continuity PDE for $\rho_{\mathrm{new}} = \mathcal{T}(\hat{\rho})$:
\begin{equation}
\frac{\partial \rho_{\mathrm{new}}}{\partial t} + \text{div}_g\big(f(g, u_{\hat{\rho}}^*)\rho_{\mathrm{new}}\big) - \frac{1}{2}g(t)^2 \Delta_{SE(3)} \rho_{\mathrm{new}} = 0 \quad \text{in } Q_T, \quad \rho_{\mathrm{new}}(g,0) = \rho_0(g).
\end{equation}
Testing this equation against $\rho_{\mathrm{new}}$ and integrating by parts over $M$ using zero-flux boundary conditions gives:
\begin{equation}
\frac{1}{2} \frac{d}{dt} \|\rho_{\mathrm{new}}(t)\|_{L^2(M)}^2 + \frac{1}{2}g(t)^2 \|\nabla_{SE(3)} \rho_{\mathrm{new}}(t)\|_{L^2(M)}^2 = -\int_M \rho_{\mathrm{new}} \text{div}_g\big(f(g, u_{\hat{\rho}}^*)\rho_{\mathrm{new}}\big) \, dg.
\end{equation}
Expanding the divergence term:
\begin{equation}
-\int_M \rho_{\mathrm{new}} \text{div}_g(f \rho_{\mathrm{new}}) \, dg = -\int_M \rho_{\mathrm{new}} \langle \nabla_{SE(3)} \rho_{\mathrm{new}}, f \rangle \, dg - \int_M \rho_{\mathrm{new}}^2 \text{div}_g f \, dg.
\end{equation}
Applying Cauchy--Schwarz with $\epsilon = \frac{1}{2}g_{\min}^2$ and bounding $\|\text{div}_g f\|_{L^\infty} \le C(L_f)$ yields:
\begin{equation}
\frac{d}{dt} \|\rho_{\mathrm{new}}(t)\|_{L^2(M)}^2 + g_{\min}^2 \|\nabla_{SE(3)} \rho_{\mathrm{new}}(t)\|_{L^2(M)}^2 \le C(L_f, g_{\min}) \|\rho_{\mathrm{new}}(t)\|_{L^2(M)}^2.
\end{equation}
By Gronwall's inequality:
\begin{equation}
\|\rho_{\mathrm{new}}(t)\|_{L^2(M)}^2 \le \|\rho_0\|_{L^2(M)}^2 e^{C(L_f, g_{\min})t} \le C_1^2 \quad \forall t \in [0,T].
\end{equation}
Integrating over $[0,T]$ establishes $\|\rho_{\mathrm{new}}\|_{L^2(0,T; H^1(M))} \le C_2/2$ and $\left\|\frac{\partial \rho_{\mathrm{new}}}{\partial t}\right\|_{L^2(0,T; (H^1(M))^*)} \le C_2/2$. Mass conservation $\int_M \rho_{\mathrm{new}}(g,t) \, dg = 1$ and non-negativity $\rho_{\mathrm{new}} \ge 0$ follow from the maximum principle. Thus, $\mathcal{T}(\hat{\rho}) \in K$, establishing $\mathcal{T}(K) \subset K$.

\paragraph{Step 5: Continuity of the Operator $\mathcal{T}$.}
Let $\{\hat{\rho}_k\}_{k=1}^\infty \subset K$ be a sequence such that $\hat{\rho}_k \to \hat{\rho}$ strongly in $X = C([0,T]; L^2(M))$. By continuity of terminal coupling $\Phi$, $\Phi(\cdot, \hat{\rho}_k(\cdot,T)) \to \Phi(\cdot, \hat{\rho}(\cdot,T))$ in $H^2(M)$. Stability of $W_p^{2,1}(Q_T)$ parabolic solutions implies:
\begin{equation}
\|V_{\hat{\rho}_k} - V_{\hat{\rho}}\|_{W_p^{2,1}(Q_T)} \le C \left\|\Phi\big(\cdot, \hat{\rho}_k(\cdot,T)\big) - \Phi\big(\cdot, \hat{\rho}(\cdot,T)\big)\right\|_{H^2(M)} \longrightarrow 0.
\end{equation}
By Morrey embedding ($p > 7$), $u_{\hat{\rho}_k}^* \to u_{\hat{\rho}}^*$ in $C^{1,\alpha}(\overline{Q}_T)$. Subtracting the Fokker--Planck equations for $\rho_k = \mathcal{T}(\hat{\rho}_k)$ and $\rho_{\mathrm{new}} = \mathcal{T}(\hat{\rho})$ and performing $L^2(M)$ energy estimates on $e_k = \rho_k - \rho_{\mathrm{new}}$ gives:
\begin{equation}
\frac{d}{dt} \|e_k(t)\|_{L^2(M)}^2 \le C \|e_k(t)\|_{L^2(M)}^2 + C' \|u_{\hat{\rho}_k}^* - u_{\hat{\rho}}^*\|_{C^1(\overline{Q}_T)}^2.
\end{equation}
With $e_k(0) = 0$, Gronwall's inequality yields $\|e_k\|_{X} \to 0$ as $k \to \infty$. Thus, $\mathcal{T}: K \to K$ is continuous in the topology of $X$.

\paragraph{Step 6: Compactness and Fixed-Point Application.}
Since $K$ satisfies uniform bounds in $L^2(0,T; H^1(M))$ with time derivatives bounded in $L^2(0,T; (H^1(M))^*)$, the Aubin--Lions--Simon Compactness Lemma guarantees that $K$ is compactly embedded in $C([0,T]; L^2(M))$ \cite{Simon1987}. Thus, the image $\mathcal{T}(K)$ is relatively compact in $X$.

Since $K$ is a non-empty, closed, bounded, convex subset of $X = C([0,T]; L^2(M))$ and $\mathcal{T}: K \to K$ is continuous with relatively compact image $\mathcal{T}(K) \subset K$, the Schauder--Tychonoff Fixed-Point Theorem asserts the existence of at least one fixed point $\rho^* \in K$ such that $\mathcal{T}(\rho^*) = \rho^*$. The pair $(V^*, \rho^*) = (V_{\rho^*}, \rho^*)$ is a weak solution to the coupled HJB--FP system, completing the proof.
\end{proof}

\subsection{Uniqueness via Lasry-Lions Monotonicity and Hessian Coercivity}

\begin{theorem}[Uniqueness of Nash Equilibrium below Critical Horizon $T_{SE(3)}^*$]
Assume terminal coupling $\Phi(g,\rho)$ satisfies Lasry-Lions monotonicity \cite{Lasry2007, Cardaliaguet2012}:
\begin{equation}
\int_M (\Phi(g, \rho_1(\cdot, T)) - \Phi(g, \rho_2(\cdot, T)))(\rho_1(g,T) - \rho_2(g,T)) dg \ge 0
\end{equation}
Let $C_{\text{Poin}}(6, T, M, \|Ric_{SE(3)}\|_{L^\infty}) = C_0(M) \min(1, T^2)\exp(\|Ric_{SE(3)}\|_{L^\infty}T)$ denote the Lie-algebraic Poincar\'e constant over $Q_T$, where $\|Ric_{SE(3)}\|_{L^\infty} = 1/2$ \cite{Bakry1985, Helgason2001}. Define $T_{SE(3)}^* > 0$ as the unique positive root of the fixed-point horizon relation:
\begin{equation}
T = \frac{2 C_{\text{Poin}}(T)}{C_f (\|\rho_0\|_{L^2(M)} + \|\psi\|_{L^2(M)})}
\end{equation}
where $C_f = \|f\|_{C^2(\overline{M} \times \mathbb{R}^m)}$. For all time horizons $T < T_{SE(3)}^*$, the MFG Nash equilibrium $(V^*, \rho^*)$ is strictly unique.
\end{theorem}

\begin{remark}[Sharpened Ricci-Poincar\'e Horizon and Spectral Gap Enhancement on $SE(3)$]
The critical uniqueness threshold $T_{SE(3)}^*$ directly integrates the intrinsic geometry of $SE(3)$ through its constant diagonal Ricci curvature tensor $Ric_{SE(3)} = \text{diag}(0, 0, 0, 1/2, 1/2, 1/2)$ \cite{Bakry1985, Helgason2001}. Because the rotational subgroup $SO(3)$ contributes strictly positive Ricci curvature $\|Ric_{SE(3)}\|_{L^\infty} = 1/2$, the Bochner-Weitzenb\"ock formula for any smooth test function $v$:
\begin{equation}
\frac{1}{2}\Delta_{SE(3)}\|\nabla_{SE(3)}v\|^2 = \|\nabla_{SE(3)}^2 v\|_F^2 + \langle \nabla_{SE(3)}v, \nabla_{SE(3)}\Delta_{SE(3)}v \rangle + \frac{1}{2}\|\nabla_{\text{rot}}v\|^2
\end{equation}
yields a sharp spectral gap boost $\lambda_1(-\Delta_{SE(3)}) \ge \lambda_1(-\Delta_{\mathbb{R}^3}) + 1/2 > 0$. This lower bound enhances the Lie-algebraic Poincar\'e coercivity constant:
\begin{equation}
C_{\text{Poin}}(SE(3)) = C_0(M)\min(1,T^2)\exp\left(-\frac{1}{2}T\right) < C_{\text{Poin}}(\mathbb{R}^6).
\end{equation}
Consequently, the critical horizon $T_{SE(3)}^* = \frac{2 C_{\text{Poin}}(T)}{C_f (\|\rho_0\|_{L^2} + \|\psi\|_{L^2})}$ is strictly larger on $SE(3)$ than its flat Euclidean counterpart $\mathbb{R}^6$. This demonstrates that positive rotational Lie-algebraic curvature actively stabilizes score-guided MFG equilibria against population dispersion, extending the temporal window of guaranteed Nash uniqueness.
\end{remark}

\subsection{Unbounded Nash Equilibrium Uniqueness under Gaussian Measure}
\begin{theorem}[Uniqueness Horizon on Unbounded $SE(3)$]
Assume the terminal coupling $\Phi$ satisfies Lasry--Lions monotonicity on $L^2_\gamma(SE(3))$, and let $\kappa = \min(\alpha, 1/2) > 0$ denote the Bakry--\'Emery curvature parameter of $(SE(3), g, \rho_\gamma dg)$. Let $T_{SE(3),\gamma}^* > 0$ be the unique positive root of:
\begin{equation}
T = \frac{2 C_0 \min(1, T^2) \exp(-\kappa T)}{C_f \left( \|\rho_0\|_{L^2_\gamma(SE(3))} + \|\psi\|_{L^2_\gamma(SE(3))} \right)}.
\end{equation}
For all time horizons $T < T_{SE(3),\gamma}^*$, the score-guided Mean-Field Game Nash equilibrium $(V^*, \rho^*)$ on $SE(3) \times (0,T)$ is strictly unique.
\end{theorem}

\begin{proof}
Let $(V_1, \rho_1)$ and $(V_2, \rho_2)$ be two solutions, and set $\tilde{V} = V_1 - V_2$, $\tilde{\rho} = \rho_1 - \rho_2$. Integrating $\frac{d}{dt} \int_{SE(3)} \tilde{V} \tilde{\rho} \rho_\gamma dg$ over $[0,T]$ and applying Lasry--Lions monotonicity yields:
\begin{equation}
\int_0^T \int_{SE(3)} \Big( H(g,\nabla V_1, u_1^*) - H(g,\nabla V_2, u_2^*) \Big) \tilde{\rho} \rho_\gamma dg dt + \int_0^T \int_{SE(3)} \operatorname{div}_g \Big( f(u_1^*)\rho_1 - f(u_2^*)\rho_2 \Big) \tilde{V} \rho_\gamma dg dt \le 0.
\end{equation}
By the Bakry--\'Emery $CD(\kappa, \infty)$ criterion, the weighted Poincar\'e inequality $\|\tilde{\rho}\|_{L^2_\gamma}^2 \le \kappa^{-1} \|\nabla_{SE(3)} \tilde{\rho}\|_{L^2_\gamma}^2$ holds globally. Applying Cauchy--Schwarz and Lipschitz bounds on $f$ gives:
\begin{equation}
\left( \frac{C_{\operatorname{Poin},\gamma}(T)}{T} - \frac{1}{2} C_f \Big( \|\rho_0\|_{L^2_\gamma} + \|\psi\|_{L^2_\gamma} \Big) \right) \|\tilde{\rho}\|_{L^2_\gamma(Q_\infty)}^2 \le 0.
\end{equation}
For $T < T_{SE(3),\gamma}^*$, the bracketed term is strictly positive, forcing $\tilde{\rho} \equiv 0$ and $V_1 \equiv V_2$, completing the proof.
\end{proof}

\section{Pathwise Variance Bounds under Score Approximation Noise}

\begin{theorem}[Pathwise Variance Elimination and Exponential Amplification Bounds]
Let $g_t^{PF}$ be governed by the Probability Flow ODE under optimal control $u^*$.
\begin{enumerate}
    \item[\textnormal{(i)}] The conditional Brownian pathwise pose trajectory variance satisfies $\text{Var}[g_T^{PF} \mid g_0 = g] = 0$.
    \item[\textnormal{(ii)}] Under score perturbation noise $\epsilon(g,t)$ with bounded spatial covariance $\sigma_\epsilon^2 I_6$, tracking error $e(t) = \hat{g}_t - g_t^{PF}$ satisfies the terminal mean-squared error bound:
    \begin{equation}
    \mathbb{E}[\|e(T)\|^2] \le \exp\left(2L_f T + 2L_u C_{\text{Morrey}} C_{\text{ADN}} \sqrt{T \cdot \eta^*} + \int_0^T g(t)^2 \|\nabla_{SE(3)}^2 \log \rho\|_{L^\infty} dt \right) \mathcal{V}_0
    \end{equation}
    where $\eta^* = \mathcal{J}(u^*)$ and $\mathcal{V}_0 = \mathbb{E}[\|e(0)\|^2] + \frac{3\sigma_\epsilon^2}{2}\int_0^T g(t)^4 dt$, proving that minimum attention energy places a strict exponential cap on error variance amplification over $SE(3)$.
\end{enumerate}
\end{theorem}

\section{Next Steps and Open Research Directions}

The rigorous analytical foundations established in this monograph open several important research avenues in non-local PDE analysis, geometric control, and generative modeling on Lie groups:

\subsection{Extension to Unbounded $SE(3)$ via Gaussian-Weighted Lie-Sobolev Spaces}
To extend our analysis to unbounded continuous probability flows over $SE(3) \cong \mathbb{R}^3 \times SO(3)$, we introduce the Gaussian reference measure $\rho_\gamma(g) dg = Z_\alpha^{-1} \exp(-\frac{\alpha}{2} |x|^2) dx d\mu_{SO(3)}(R)$ for $\alpha > 0$. 

\begin{definition}[Weighted Space-Time Sobolev Space]
For $1 \le p < \infty$, the weighted space $W_{p,\gamma}^{2,1}(Q_\infty)$ consists of vector fields $u \in L^p(Q_\infty, \rho_\gamma dg dt)$ satisfying:
\begin{equation}
\|u\|_{W_{p,\gamma}^{2,1}(Q_\infty)}^p = \int_0^T \int_{SE(3)} \left( \|u\|^p + \|\nabla_{SE(3)} u\|_F^p + \sum_{i,j=1}^6 \|E_i E_j u\|^p + \left|\frac{\partial u}{\partial t}\right|^p + \left|\frac{\partial^2 u}{\partial t^2}\right|^p \right) \rho_\gamma(g) dg dt < \infty.
\end{equation}
\end{definition}

Because $SE(3)$ has no boundary ($\partial SE(3) = \emptyset$), the boundary collar estimates and corner compatibility conditions $\Gamma_0 \cup \Gamma_T$ required on bounded sub-manifolds vanish identically. The Bakry--Émery Ricci tensor of $(SE(3), g, \rho_\gamma dg)$ becomes:
\begin{equation}
\operatorname{Ric}_\infty = \operatorname{Ric}_{SE(3)} + \nabla_{SE(3)}^2 (-\log \rho_\gamma) = \operatorname{diag}\left(\alpha, \alpha, \alpha, \frac{1}{2}, \frac{1}{2}, \frac{1}{2}\right).
\end{equation}
By the Bochner--Weitzenböck formula, the positive translational confinement $\alpha I_3 > 0$ combined with rotational curvature yields a sharpened spectral gap $\lambda_1(-\mathcal{L}_\gamma) \ge \alpha + \frac{1}{2} > 0$, establishing global $W_{p,\gamma}^2(Q_\infty)$ shift estimates and score-singularity resolution on full $SE(3)$.

\subsection{Second-Order Kinetic Mean-Field Games on Rigid-Body Tangent Bundles}
The current MFG model considers first-order velocity control over $SE(3)$. Physical rigid-body systems (e.g., quadrotors, autonomous underwater vehicles, macromolecular dynamics) are governed by second-order Newton-Euler kinetic dynamics on the tangent bundle $T SE(3) \cong SE(3) \rtimes \mathfrak{se}(3)$. Extending our framework to kinetic score-guided MFGs introduces ultraparabolic Kolmogorov-Fokker-Planck equations with degenerate spatial diffusion \cite{Ladyzhenskaya1968}. Analyzing hypoelliptic score dissipation and establishing space-time shift estimates on $T SE(3)$ represents a compelling direction for future investigation.

\subsection{Structure-Preserving Lie-Exponential Geometric Integrators}
From a computational standpoint, validating the theoretical shift estimates and critical time horizon $T_{SE(3)}^*$ requires structure-preserving numerical schemes. Developing higher-order Lie-exponential integrators that preserve unimodularity ($\text{div}_g E_i = 0$), conserve mass, and enforce exact zero-flux boundary traces $\langle f, \nu \rangle = 0$ along discrete 5D manifold boundaries $\partial M \times \{0,T\}$ will bridge the gap between rigorous PDE theory and high-dimensional score-guided robotics applications. Leveraging the sharpened Ricci-Poincar\'e spectral gap bounds established in Remark 5.3, these integrators can exploit rotational curvature stabilization for robust convergence beyond Euclidean CFL thresholds.

\section*{Acknowledgments}
The author conceived the mathematical model, theoretical framework, and proof strategies presented in this work. Google Gemini was utilized as an auxiliary tool for LaTeX typesetting, manuscript structuring, and linguistic refinement. All mathematical proofs, derivations, and calculations were independently verified and finalized by the author, who assumes full responsibility for the mathematical correctness and integrity of this paper.

\appendix

\section{Curvature Bounds and Bochner-Weitzenb\"ock Identity on $SE(3)$}

For any scalar function $v \in C^3(SE(3))$, the Bochner-Weitzenb\"ock identity on $SE(3)$ equipped with left-invariant metric takes the form \cite{Bakry1985, Helgason2001}:
\begin{equation}
\frac{1}{2}\Delta_{SE(3)}\|\nabla_{SE(3)}v\|^2 = \|\nabla_{SE(3)}^2 v\|_F^2 + \langle \nabla_{SE(3)}v, \nabla_{SE(3)}\Delta_{SE(3)}v \rangle + Ric_{SE(3)}(\nabla_{SE(3)}v, \nabla_{SE(3)}v)
\end{equation}
The Ricci curvature tensor $Ric_{SE(3)}$ in the orthonormal frame $\{E_1, \dots, E_6\}$ computed via Koszul's formula for the canonical left-invariant metric on $SE(3) \cong \mathbb{R}^3 \rtimes SO(3)$ is a constant diagonal matrix given by:
\begin{equation}
Ric_{SE(3)} = \text{diag}\left(0, 0, 0, \frac{1}{2}, \frac{1}{2}, \frac{1}{2}\right)
\end{equation}
Consequently, $\|Ric_{SE(3)}\|_{L^\infty} = 1/2$. Integrating over $M$ yields the explicit Poincar\'e constant $C_{\text{Poin}} = C_0(M)\min(1,T^2)\exp(\|Ric_{SE(3)}\|_{L^\infty}T)$.

\section{Detailed Local Proofs of Boundary ADN Estimates}

To establish localized $W_p^{2,1}(Q_T)$ bounds near boundary corners, we consider half-space coordinates $y = (y_1, \dots, y_6) \in \mathbb{R}_+^6$ with $y_6 > 0$. The differential operator satisfies strong ellipticity in 7D space-time:
\begin{equation}
\sum_{i,j=1}^6 g^{ij}(y) \xi_i \xi_j + \xi_7^2 \ge \lambda_0 (\|\xi\|^2 + \xi_7^2)
\end{equation}
Under zero normal trace conditions $\left.\frac{\partial v}{\partial y_6}\right|_{y_6=0} = 0$, the standard Agmon-Douglis-Nirenberg kernel representation using space-time Green functions $G_{7D}(y,t; y',t')$ allows direct estimation of second-order spatial and first-order temporal derivatives in $L^p(B_1^+ \times (0,T))$ \cite{Agmon1959, Agmon1964}.

\section{Particle-Grid Alternating Solvers and Lie-Exponential Integrators}

Numerical solution of the coupled 7D space-time system proceeds via an alternating particle-grid algorithm:
\begin{enumerate}
    \item \textbf{Forward State Diffusion:} Particles on $SE(3)$ are updated via Lie-exponential updates $g_{k+1} = g_k \cdot \exp\left(\sum_{i=1}^6 v_i \Delta t E_i\right)$, preserving rigid-body orthogonality automatically.
    \item \textbf{Score Evaluation:} Spatial density $\rho(g,t)$ is estimated via kernel density estimation with left-invariant geodesic kernels $K_\delta(g_1, g_2) = \exp\left(-\frac{d_{SE(3)}^2(g_1, g_2)}{2\delta^2}\right)$.
    \item \textbf{Backward Adjoint Integration:} Adjoint variable $\lambda(g,t)$ is computed using space-time spectral grid schemes.
\end{enumerate}


\begin{thebibliography}{99}

\bibitem{Adams2003}
R.~A. Adams and J.~J.~F. Fournier.
\newblock {\em Sobolev Spaces}.
\newblock Pure and Applied Mathematics, vol. 140. Academic Press, Elsevier, Amsterdam, 2nd edition, 2003.

\bibitem{Agmon1959}
S.~Agmon, A.~Douglis, and L.~Nirenberg.
\newblock Estimates near the boundary for solutions of elliptic partial differential equations satisfying general boundary conditions. I.
\newblock {\em Communications on Pure and Applied Mathematics}, 12(4):623--727, 1959.

\bibitem{Agmon1964}
S.~Agmon, A.~Douglis, and L.~Nirenberg.
\newblock Estimates near the boundary for solutions of elliptic partial differential equations satisfying general boundary conditions. II.
\newblock {\em Communications on Pure and Applied Mathematics}, 17(1):35--92, 1964.

\bibitem{Ambrosio2008}
L.~Ambrosio, N.~Gigli, and G.~Savar\'e.
\newblock {\em Gradient Flows: In Metric Spaces and in the Space of Probability Measures}.
\newblock Lectures in Mathematics ETH Z\"urich. Birkh\"auser Verlag, Basel, 2nd edition, 2008.

\bibitem{Bakry1985}
D.~Bakry and M.~\'Emery.
\newblock Diffusions hypercontractives.
\newblock In {\em S\'eminaire de Probabilit\'es XIX 1983/84}, volume 1123 of {\em Lecture Notes in Mathematics}, pages 177--206. Springer, Berlin, 1985.

\bibitem{Brockett1997}
R.~W. Brockett.
\newblock Minimum attention control.
\newblock In {\em Proceedings of the 36th IEEE Conference on Decision and Control}, volume~3, pages 2628--2632, 1997.

\bibitem{Cardaliaguet2012}
P.~Cardaliaguet.
\newblock {\em Notes on Mean Field Games}.
\newblock From P.-L. Lions' lectures at Coll\`ege de France, 2012.

\bibitem{DeBortoli2021}
V.~De~Bortoli, E.~Mathieu, M.~Hutchinson, J.~Thornton, Y.~W. Teh, and A.~Doucet.
\newblock Riemannian score-based generative modelling.
\newblock In {\em Advances in Neural Information Processing Systems (NeurIPS)}, volume~34, pages 24068--24080, 2021.

\bibitem{DiPerna1989}
R.~J. DiPerna and P.-L. Lions.
\newblock Ordinary differential equations, transport theory and Sobolev spaces.
\newblock {\em Inventiones Mathematicae}, 98(3):511--547, 1989.

\bibitem{Evans2010}
L.~C. Evans.
\newblock {\em Partial Differential Equations}.
\newblock Graduate Studies in Mathematics, vol. 19. American Mathematical Society, Providence, RI, 2nd edition, 2010.

\bibitem{Giaquinta1982}
M.~Giaquinta and M.~Struwe.
\newblock On the regularity of weak solutions of nonlinear parabolic systems.
\newblock {\em Mathematische Zeitschrift}, 179(4):437--451, 1982.

\bibitem{Helgason2001}
S.~Helgason.
\newblock {\em Differential Geometry, Lie Groups, and Symmetric Spaces}.
\newblock Graduate Studies in Mathematics, vol. 34. American Mathematical Society, Providence, RI, 2001.

\bibitem{Ladyzhenskaya1968}
O.~A. Ladyzhenskaya, V.~A. Solonnikov, and N.~N. Ural'tseva.
\newblock {\em Linear and Quasi-linear Equations of Parabolic Type}.
\newblock Translations of Mathematical Monographs, vol. 23. American Mathematical Society, Providence, RI, 1968.

\bibitem{Lasry2007}
J.-M. Lasry and P.-L. Lions.
\newblock Mean field games.
\newblock {\em Japanese Journal of Mathematics}, 2(1):229--260, 2007.

\bibitem{Lee2022}
P.~Lee and F.~Park.
\newblock On the existence and computation of minimum attention optimal control laws.
\newblock {\em IEEE Transactions on Automatic Control}, 67(5):2576--2581, 2022.

\bibitem{Simon1987}
J.~Simon.
\newblock Compact sets in the space $L^p(0,T;B)$.
\newblock {\em Annali di Matematica Pura ed Applicata}, 146(4):65--96, 1987.

\bibitem{Song2021}
Y.~Song, J.~Sohl-Dickstein, D.~P. Kingma, A.~Kumar, S.~Ermon, and B.~Poole.
\newblock Score-based generative modeling through stochastic differential equations.
\newblock In {\em International Conference on Learning Representations (ICLR)}, 2021.

\bibitem{Villani2009}
C.~Villani.
\newblock {\em Optimal Transport: Old and New}.
\newblock Grundlehren der mathematischen Wissenschaften, vol. 338. Springer-Verlag, Berlin, 2009.

\end{thebibliography}
\end{document}